\documentclass[
 amsmath,amssymb,
 preprint
]{revtex4-1}

\usepackage{float}
\usepackage{graphicx}
\usepackage{dcolumn}
\usepackage{bm}
\usepackage{subfig}
\usepackage{amsthm, amssymb, amsmath}

	\newtheorem{theorem}{Theorem}[section]
	\newtheorem{definition}[theorem]{Definition}
	\newtheorem{proposition}[theorem]{Proposition}
	\newtheorem{lemma}[theorem]{Lemma}
	\newtheorem{remark}[theorem]{Remark}
	\newtheorem{sublemma}[theorem]{Sublemma}
	\newtheorem{corollary}[theorem]{Corollary}
	\newtheorem{assumption}[theorem]{Assumption}
	\newtheorem{notationalrem}[theorem]{Notational Remark}
	\newtheorem{conjecture}[theorem]{Conjecture}   
	\newtheorem{tools}[subsection]{$\negsp\negsp$}
	\newcommand\asm[1]{ \begin{assumption}\label{#1} }
	\newcommand\easm{ \end{assumption} }
	\newcommand\dfn[1]{ \begin{definition}\label{#1} }
	\newcommand\dfntwo[2]{ \begin{definition}[#2]\label{#1} }
	\newcommand\edfn{ \end{definition} }
	\newcommand\rem[1]{ \begin{remark}\label{#1} \rm}
	\newcommand\remtwo[2]{ \begin{remark}[#2]\label{#1} \rm}
	\newcommand\erem{ \end{remark} }
	\newcommand\thm[1]{ \begin{theorem}\label{#1}}
	\newcommand\thmtwo[2]{ \begin{theorem}[#2]\label{#1}}
	\newcommand\ethm{ \end{theorem} }
	\newcommand\pro[1]{ \begin{proposition}\label{#1}}
	\newcommand\protwo[2]{ \begin{proposition}[#2]\label{#1}}
	\newcommand\epro{ \end{proposition} }
	\newcommand\lem[1]{ \begin{lemma}\label{#1}}
	\newcommand\lemtwo[2]{ \begin{lemma}[#2]\label{#1}}
	\newcommand\elem{ \end{lemma} }
	\newcommand\sublem[1]{ \begin{sublemma}\label{#1}}
	\newcommand\sublemtwo[2]{ \begin{sublemma}[#2]\label{#1}}
	\newcommand\esublem{ \end{sublemma} }
	\newcommand\cor[1]{ \begin{corollary}\label{#1}}
	\newcommand\cortwo[2]{ \begin{corollary}[#2]\label{#1}}
	\newcommand\ecor{ \end{corollary} }
	\newcommand\notrem[1]{ \begin{notationalrem}\label{#1} \sl}
	\newcommand\enotrem{ \end{notationalrem} }

	\newtheorem{algorithm}[theorem]{Algorithm}

	\newcommand\avg[1]{{ \int_{\mathbb{T}} #1 d\theta}}

	\newcommand\equ[1]{{\rm (\ref{#1})}}
	\newcommand\beq[1]{ \begin{equation}\label{#1} }
	\newcommand{\eeq}{ \end{equation} }

	\newcommand\beqa[1]{ \begin{eqnarray} \label{#1}}
	\newcommand{\eeqa}{ \end{eqnarray} }
	\newcommand{\beqano}{ \begin{eqnarray*} }
	\newcommand{\eeqano}{ \end{eqnarray*} }

	\renewcommand{\Im}{{\rm \, Im\,}}

	\newcommand{\negsp}{\hspace{-.04truecm}}

	\newcommand{\e }{ {\varepsilon}   }

	\font\teneufm=eufm10
	\font\seveneufm=eufm7
	\font\fiveeufm=eufm5
	\newfam\eufmfam
	\textfont\eufmfam=\teneufm
	\scriptfont\eufmfam=\seveneufm
	\scriptscriptfont\eufmfam=\fiveeufm

\begin{document}


\title{\bf Computation of Domains of Analyticity for the dissipative standard
map in the limit of small dissipation}

\author{Adri\'an  P. Bustamante}
\email{apb7@math.gatech.edu}
\affiliation{School of Mathematics \\
Georgia Institute of Technology\\
686 Cherry St. Atlanta, GA, 30332-0160, USA}

\author{Renato C. Calleja}
\email{calleja@mym.iimas.unam.mx}
\affiliation{Department of Mathematics and Mechanics\\
IIMAS, National Autonomous University of Mexico (UNAM)\\
Apdo. Postal 20-126, C.P. 01000, Mexico, D.F. MEXICO}

\date{\today}

\begin{abstract}
\noindent
Conformally symplectic systems include mechanical
systems with a friction proportional to the velocity.
Geometrically, these systems transform a symplectic form into 
a multiple of itself making the systems dissipative or expanding.
In the present work we consider the limit of small dissipation. 
The example we study is a family of conformally symplectic standard
maps of the cylinder
for which the conformal factor, $b(\e)$,
is a function of a small complex parameter, $\e$.

We assume that for $\e=0$ the map preserves the symplectic form and
the 
dependence on $\e$ is cubic, i.e., $b(\e) = 1 - \e^3$.
We compute perturbative expansions formally in $\e$ and use them
to estimate the shape of the domains of analyticity of invariant circles
as functions of $\e$.
We also give evidence that the functions might belong to a Gevrey
class at $\e = 0$.
We also perform numerical continuation of the solutions
as they pass through the boundary of the domain to illustrate
that the monodromy of the solutions is trivial. 
The numerical computations we perform support conjectures
on the shape of the domains of analyticity. 

\end{abstract}


\maketitle

\section{Introduction}

We study the limit of small dissipation/expansion of a family 
of conformally symplectic standard maps. 
In particular, we approximate the shape of domains of analyticity of invariant
circles of a family of conformally symplectic standard maps
of the cylinder, $\mathcal M = \mathbb S^1 \times \mathbb R$, depending 
on a small parameter, $\e$, that vanishes as the conformal 
factor tends to one.

It was noted in \cite{Cal-Cel-Lla-16}
that the small divisors depend on the complex
parameter $\e$ and give rise to regions where the functions parameterizing
the circles cannot be analytic with respect to $\e$
but miss by very little.
A conjecture in \cite{Cal-Cel-Lla-16} states that the tori
are analytic in a domain in the complex $\e$ plane that is obtained by
taking from a ball centered at zero, a sequence of small balls with
centers along smooth curves passing through the origin.
The radii of the excluded balls decreases faster that any power of the
distance of the centers of the balls to the origin. In fact, it was rigorously proved in
\cite{Cal-Cel-Lla-16} that this domain is a lower bound.
The main objective of the present work is to
illustrate the results in one example, provide numerical evidence and indications
of new results.
Our computations indicate that there are singularities which
cluster around several points at which one does not expect the 
functions to be analytic. The singularities in the complex $\e$ plane,
cluster inside balls whose radii decrease at the ratios
predicted by the conjecture.

A common method to compute invariant circles of a map of the cylinder
$f_\e : \mathcal M \to \mathcal M$, is by computing a parameterization
$K_\e : \mathbb S^1 \to \mathcal M$
of the invariant circle which
satisfies an invariance equation. The invariance equation is
\[f_\e \circ K_\e = K_\e \circ T_\omega\]
with $T_\omega(\theta) = (\theta + \omega)$.
The invariance equation states that the dynamics on the
invariant circle are conjugated to a rigid rotation
of the circle by an irrational number $\omega$. 
The parameterization function $K_\e$ can be written in terms of a periodic
function, $u_\e : \mathbb S^1 \to \mathbb R$, as in equation \eqref{parameterization}.
The method we use to find the singularities is to approximate the conjugacy function
$u_\e(\theta)$ by means of a Lindstedt series expansion in $\e$.
The Lindstedt method produces polynomials in $\e$ of high order,
\[u^{\leq N}_\varepsilon(\theta)=\sum_{n=0}^N u_n(\theta)\varepsilon^n\]
with $N\approx 10^3$, 
whose coefficients $u_n: \mathbb{S}^1 \to \mathbb C$ are periodic
functions. 
We then use the Lindstedt series of the conjugacies to obtain Pad\'e rational functions whose poles are expected approximate the poles of the original
function $u_\e$. Pad\'e extrapolation methods of Lindstedt series have been widely
used by several authors
\cite{Ber-Fal-Gen-01, Ber-Gen-01, Cel-Fal-02, Ber-Gen-04, Lla-Tom-95}
in the simplectic case.
Since the Pad\'e extrapolation method is based on
approximating an analytic function with a rational function,
the computation of poles is done by approximating the roots of the denominator
of the Pad\'e function. The denominator is a polynomial that can be of
very high degree,
and computing its roots depends heavily on numerical precision.
Since the computations are very
sensitive to precision, we perform them using $\approx 10^3$ digits which allows us
to compute singularities for values of $\e$ 
that are at a distance $\approx 0.3$ from $\e = 0$
in the complex plane.
We expect that higher precision together with higher order degree series,
would allow us to compute poles that are closer
to the origin. However, the method already allows us to have an approximation
of the boundary of the domain
in regions that are contained in the small balls that were predicted by the
conjecture in \cite{Cal-Cel-Lla-16}, even when the singularities are not
very close to $\e = 0$. For this reason it is very hard to notice that the functions that
we are computing are not analytic. We also find conjectures on the rate of growth of the terms of the Lindstedt series.

We note that the shapes of the domains that we present here are 
remarkably different from what one sees in the symplectic case,
see \cite{Cel-Fal-02, Ber-Gen-04, Lla-Tom-95, Cal-Lla-10}.
This is
partly due to the fact that in the symplectic case
or in the dissipative case, the small divisors
do not depend on the conformal factor $b(\e)$ which in our case is
a function of $\e$.

Some explorations of the shape of the analyticity domains in
the dissipative standard map have been performed using the parameterization
method in \cite{Cal-Cel-10}, that is very similar
to the one described in section \ref{newton}.
In \cite{Cal-Cel-10}, it is noticed that the breakdown of invariant
tori in the conservative and the dissipative case are similar
when the conformal factor $b$ is a constant,
\cite{Cal-Lla-10b}. A different behavior in the breakdown of invariant tori involving bundle 
collapse is observed in the dissipative standard map in \cite{Cal-Fig-12}.
Explorations of the shape of the domains of analyticity
in $\e$ in the conservative case
with the use of the parameterization method appear in 
\cite{Cal-Lla-09, Cal-Lla-10}.

\section{Preliminaries}
We consider the dissipative standard map
defined on the cylinder $\mathcal M = \mathbb S^1 \times \mathbb R$ given
by $f_\e(x_n, y_n) = (x_{n+1}, y_{n+1})$ and 
\beqa{smdiss} \label{std-dis-ham}
y_{n+1}&=&b_\e y_n+ c_\e +{\varepsilon}\ V'(x_n) \nonumber\\  
x_{n+1}&=&x_n+y_{n+1}\ , 
\eeqa
where $y_n\in\mathbb{R}$, $x_n\in\mathbb{S}^1$, $\varepsilon\in \mathbb{R}$, and $V'(x)=\frac{1}{2\pi}\sin(2\pi x)$ is an
analytic, periodic function. Here we consider the case
when the dissipative parameter, $b_\e$, is given by $b_\e=b(\e)=1-\varepsilon^3$, and
the drift parameter $c=c(\varepsilon)$ is a function that depends on the small parameter
$\e$. The dissipative parameter $b_\e$ coincides with the Jacobian of the function.
We note that the Jacobian is the rate of dissipation/expansion of the map
\eqref{std-dis-ham}, this rate will be dependent
of the parameter $\e$. In particular, the case $\e = 0$ coincides with the
zero dissipation limit. 

In fact, it is discussed in \cite{Cal-Cel-Lla-16} that \eqref{std-dis-ham}
is conformally symplectic.
If $\Omega=dy \wedge dx$ is the standard symplectic form of the cylinder, the map
$f_\e$ satisfies that
\begin{equation}\label{conf-symp}
  f^*\Omega = b_\e \Omega .
\end{equation}
For certain values of $c_\e$, we know that maps of the from \eqref{conf-symp}
have analytic
invariant circles corresponding to quasi-periodic orbits
with Diophantine rotation numbers, $\omega$.
The Lindstedt series analysis in Section \ref{lindstedt} determines
that one condition for the
the mapping \eqref{std-dis-ham} to admit an invariant circle is that
$c_\e = \omega \e^3 + \mathcal O(\e^4)$. In the following, we discuss the properties
that the rotation number should satisfy so that one can have quasi-periodic
orbits parameterized by a function. 

\subsection{Quasi-periodic orbits}

We consider a frequency $\omega$ that satisfies the Diophantine condition,
\begin{equation}\label{dioph-cond}
  |\omega q - p|\geq \nu |q|^{-\tau}, \quad p \in \mathbb Z, \quad q \in \mathbb Z\setminus \{0\}
\end{equation}
where $\nu \in \mathbb R^+$ and $\tau \in \mathbb R$ with $\tau \geq 1$.

Quasi periodic orbits of the dissipative standard map \eqref{std-dis-ham} are found
using a parametric representation of the variable $x_n \in \mathcal{S}^1$ as
\begin{equation}\label{conj-rep}
  x_n = \theta_n + u_\e(\theta_n), \quad \theta\in \mathbb S^1,
\end{equation}
where $u_\e:\mathbb S^1 \to \mathbb R$ is a $1$-periodic function. We assume that
the variable $\theta_n$ varies linearly as $\theta_{n+1} = \theta_{n} + \omega$ where
$\omega$ is the rotation frequency.

It follows form equation \eqref{std-dis-ham} that
\begin{equation}\label{std-dis-lag}
  x_{n+1} -(1+b_\e)x_n +b_\e x_{n-1} +(1-b_\e)\omega - c_\e +\varepsilon V'(x_n) = 0
\end{equation}
We look for quasi periodic solutions by finding $u_\e$
and $c_\e = c(\e)$ such that
\begin{equation}\label{func-eqn}
  E_{c_\e}[u_\e] \equiv u_\varepsilon(\theta +\omega) -(1+b_\e)u_\varepsilon(\theta)
  +b_\e u_\varepsilon(\theta -\omega) +(1-b_\e)\omega -c_\e
  + \varepsilon V'(\theta +u_\varepsilon(\theta)) = 0.
\end{equation}
We remark that the nature of the two unknowns is different since $u_\e(\theta)$ is a
smooth complex $1$-periodic function of $\theta \in \mathbb S^1$  depending on the complex
parameter $\e$ and $c_\e$ is a complex number depending on $\e$.
The conjecture in \cite{Cal-Cel-Lla-16}, states that $\e$ is a complex parameter
whose range lays in a complex domain that is
obtained by taking out from a neighborhood of $\e = 0$,
points inside balls with centers along smooth curves passing though
the origin.
In \cite{Cal-Cel-Lla-16} there is also a rigorous lower bound close to the domain
desribed in the conjecture.

It is clear that once we find a pair $(u_\e, c_\e)$ satisfying \eqref{func-eqn}, we can
recover the embedding of the quasi-periodic orbit by the parameterization
$K_\e : \mathbb S^1 \to \mathcal M$,
\begin{equation}\label{parameterization}
K_\e(\theta) = \left(\begin{array}{c}
\theta + u_\e (\theta) \\
\omega + u_\e(\theta) - u_\e(\theta - \omega)\\
\end{array}\right).
\end{equation}

\section{Methods for computing solutions}\label{methods}

We will use two different methods for finding the solution pair $(u_\e, c_\e)$
of \eqref{func-eqn}. The first method is based on a Lindstedt series
approximation of the solutions written as formal power series 
of the small parameter $\e$. In our case the small parameter $\e$
will account both for the 
size of the perturbation and the distance of the conformal factor
to the symplectic case.
This method produces approximate solutions in the sense that if 
\begin{equation}\label{form-pow-ser}
  u^{\leq N}_\varepsilon(\theta)=\sum_{k=0}^N u_k(\theta)\varepsilon^k\quad
  \mbox{and}\quad
  c^{\leq N}(\varepsilon) =\sum_{k=0}^N c_k\varepsilon^k
\end{equation}
are polynomials in $\e$, we say that \eqref{form-pow-ser} is an approximate solution
of order $N$ whenever $\|E_{c^{\leq N}(\varepsilon)}[ u^{\leq N}_\varepsilon]\| = O(\e^{N+1})$, where
$E$ is the functional defined in \eqref{func-eqn} and $\|\cdot \|$ is the supremum norm
over all $\theta \in \mathbb S^1$. The Lindstedt series method that we describe in section
\ref{lindstedt} provides a way to construct an approximate solution of any given
order $N\in \mathbb N$.

In section \ref{newton}, we include an algorithm to find the solution
$(u_\e, c_\e)$ by means of a Newton method. The method starts form an approximate solution pair 
$(u_a, c_a)$ so that the norm of $E_{c_a}[u_a]$ is small and provides
a correction $(v,\delta)$ by imposing that the new solution
$(u_a+v, c_a+\delta)$ satisfies the functional equation
$E_{c_a+\delta}[u_a+v]$ up to first order in $(v,\delta)$. This method
can be shown to converge using scales of Banach spaces.

\subsection{Lindstedt Series}\label{lindstedt}

The Lindstedt series method consists of performing a formal series expansion in
a small parameter $\e$. According to \eqref{func-eqn},  and the fact that $b(\e)=1-\e^3$,
we look for a solution, $(u_\e, c_\e)$, of
\begin{equation}\label{ecu_u}
u_\varepsilon(\theta+\omega) -(2-\varepsilon^3)u_\varepsilon(\theta) +(1-\varepsilon^3)u_\varepsilon(\theta-\omega) +\varepsilon^3\omega - c(\varepsilon) =-\varepsilon V'(\theta +u_\varepsilon(\theta)) 
\end{equation}
as a power series expansion. That is, 
we look for solutions 
\begin{equation}\label{form-pow-ser_u}
  u_\varepsilon(\theta)=\sum_{k=0}^\infty u_k(\theta)\varepsilon^k\quad
  \mbox{and}\quad
  c(\varepsilon) =\sum_{k=0}^\infty c_k\varepsilon^k,
\end{equation}
where each $u_n$ is a function from $\mathbb{S}^1$ to $\mathbb C$ and each $c_n\in \mathbb C$.
This solution can be computed by equating powers of $\e$ in \eqref{ecu_u}. Taking the Taylor
expansion at $\varepsilon=0$ 
\begin{equation}\label{Taylor_v}
-\varepsilon V'(\theta+u_\varepsilon(\theta)) =\sum_{k=1}^\infty S_k(\theta)\varepsilon^k
\end{equation}
and substituting \eqref{form-pow-ser_u} into \eqref{ecu_u}, we have that 
\begin{equation} \label{ecu_u_ser}\sum_{k=0}^\infty u_k(\theta+\omega)\e^k -(2-\e^3)\sum_{k=0}^\infty u_k(\theta) +(1-\e^3)\sum_{k=0}^\infty u_k(\theta-\omega) +\e^3\omega - \sum_{k=0}^\infty c_k\e^k =\sum_{k=1}^\infty S_k(\theta)\e^k.
\end{equation}

\begin{remark}\label{Rem-S_k} When $V'(\theta)=\sin(2\pi\theta)$, or a trigonometric polynomial, the $S_k(\theta)$'s can be
  computed very efficiently in terms of the $u_i(\theta)$'s. Following \cite{KNUTH, Fal-Lla-92} and denoting $\mathcal S(\theta,\e)= \sin(2\pi(\theta +u_\e(\theta)))$,
  $\mathcal C(\theta,\e)= \cos(2\pi(\theta+ u_\e(\theta)))$, the coefficients of the series expansions
  $\displaystyle{\mathcal S(\theta,\e)= \sum_{k=0}^\infty \mathcal S_k(\theta)\e^k}$ and $\displaystyle{\mathcal C(\theta,\e)
    = \sum_{k=0}^\infty \mathcal C_k(\theta)\e^k}$ are given by the
  following recurrence relations,
  \begin{eqnarray}
    (N+1)\mathcal S_{N+1}(\theta)=& 2\pi\sum_{m=0}^N \mathcal C_{N-m}(m+1)u_{m+1}(\theta) \\
    (N+1)\mathcal C_{N+1}(\theta)=&-2\pi\sum_{m=0}^N \mathcal S_{N-m}(m+1)u_{m+1}(\theta).\nonumber
  \end{eqnarray}
  Thus $S_k(\theta)=-\mathcal S_{k-1}(\theta)$ and $S_0\equiv 0$, by \eqref{Taylor_v}.
\end{remark}

Defining the operator
\begin{equation}\label{Lomega}
  L_\omega\varphi(\theta):= \varphi(\theta +\omega) -2\varphi(\theta) +\varphi(\theta -\omega)
\end{equation}
equation \eqref{ecu_u_ser} can be rewritten as    
\begin{align}
\sum_{k=1}^\infty S_k(\theta)\varepsilon^k &=\sum_{k=0}^2 \left(L_\omega u_k(\theta)-c_k\right)\varepsilon^k +\left(L_\omega u_3(\theta)-c_3 +u_0(\theta) -u_0(\theta -\omega) +\omega\right)\varepsilon^3 \nonumber\\ 
 &\quad +\sum_{k=4}^\infty \left(L_\omega u_k(\theta) -c_k +u_{k-3}(\theta) -u_{k-3}(\theta-\omega)\right)\varepsilon^k, \label{ecu_pot_lindst}
\end{align}

Some properties of the operator $L_\omega$ are summarized in the following Lemma.
See \cite{Lla-01} for details about the proof.

\lem. \label{L_w_sol} Let $\eta:\mathbb{S}^1\rightarrow\mathbb{S}^1$ a continuous function such
that $\displaystyle{\int_0^1\eta(\theta)d\theta=0}$. If $\omega$ is Diophantine
as in \eqref{dioph-cond}, then there
exists a solution, $\varphi(\theta)$, to the equation
\begin{equation}\label{solveLomega}
  L_\omega\varphi(\theta)=\eta(\theta)
\end{equation}
such
that $\displaystyle{\int_0^1\varphi(\theta)d\theta=0}$. In fact, the solution is given by
$$\varphi(\theta) =\sum_{\ell \in\mathbb{Z}\backslash \{0\}}\frac{\hat{\eta}_\ell}{2 (\cos(2\pi \ell\omega)-1)}
e^{2\pi i \ell\theta},$$ where $\hat{\eta}_\ell$ are the Fourier coefficients of $\eta(\theta)$.
\elem


The Lindstedt process is as follows. Matching the coefficients of
the same order in \eqref{ecu_pot_lindst} we
obtain the following relations to different orders of $\e$.
The zero-th order term tells us that the coefficients at order zero in $\e$ have to be
trivial. The equations are
\begin{equation}
L_\omega u_0(\theta)-c_0=0. \label{ec_u0}
\end{equation}
Choosing $c_0= 0$, then $u_0\equiv 0$ is the solution given by Lemma \ref{L_w_sol}.
This construction is analogous to the zero-th order term in the symplectic case.

\begin{remark} \label{Rem-non-dis}
  This method has been used in \cite{Fal-Lla-92, Lla-01,
      Ber-Fal-Gen-01, Ber-Gen-01, Cel-Fal-02, Ber-Gen-04, Lla-Tom-95} for the symplectic
  case. That is, making the same process for the standard map,
  $(x_{n+1}, y_{n+1})=(x_n + y_{n+1}, y_n +\e V'(x_n))$, gives the following
  equation to all orders $\e^n$, 
\begin{equation}
L_\omega u_k(\theta) =S_k(\theta)\qquad k\geq 0. 
\end{equation}
Moreover, $\displaystyle{\int_0^1 S_k(\theta)d\theta =0}$ for all $k\geq 0$.
This is a consequence of the symplectic structure and
the fact that $S_k(\theta)$ depends on the previously
computed
$u_0(\theta)$, $u_1(\theta)$, $\dots$,
$u_{k-1}(\theta)$, $S_0(\theta)$, $S_1(\theta)$, $\dots$,  $S_{k-1}(\theta)$
(see Remark \ref{Rem-S_k}).  
\end{remark}

The first and second orders in $\e$ are also analogous to the symplectic case.
For this reason the first two coefficients of $c_\e$ will be trivial.
\begin{equation}
L_\omega u_k(\theta)-c_k= S_k(\theta), \qquad k= 1,2.  \label{ec_u1}
\end{equation}
Choosing $c_1=0=c_2$ the equations are reduced to the non dissipative case and, by Remark
\ref{Rem-non-dis}, the right hand side has zero average. Therefore, we can find solutions
$u_1(\theta)$, $u_2(\theta)$.

The third order in $\e$ is the
first one that is different from the conservative case.
\begin{equation}
L_\omega u_3(\theta) -c_3 +\omega =S_3(\theta).\label{ec-u3}
\end{equation}
Here we notice that the drift parameter starts playing
a r\^ole in the existence of invariant tori.
Taking $\displaystyle{c_3 =\omega}$, equation \eqref{ec-u3} becomes the same equation as in the symplectic case. Since $S_3(\theta)$ has zero average we find
$u_3(\theta)$.

The equations for orders higher that $4$ are remarkably different since
we have a counter term coming from the previously computed orders. Namely,
\begin{equation}
  L_\omega u_k(\theta) = S_k(\theta) -u_{k-3}(\theta) +u_{k-3}(\theta-\omega)
  +c_k,\qquad k\geq 4. \label{ec-u_n}
\end{equation}
Notice that, by construction, $\displaystyle{\int_0^1u_{k-3}(\theta-\omega)d\theta
  = \int_0^1u_{k-3}(\theta)d\theta=0}$ (see Lemma \ref{L_w_sol}). Now, taking  
\begin{equation}
c_k =-\int_0^1 S_k(\theta)d\theta,\label{mu-n}
\end{equation}
we can find $u_k(\theta)$ solving \eqref{ec-u_n} for all $k\geq 4$.

We have proved the following proposition which is a particular case of part A)
of Theorem 12 in \cite{Cal-Cel-Lla-16}.

\begin{proposition}\label{prop:existenceLind}
  For any $N\in \mathbb{N}$, the procedure presented
  above allows to find an approximate solution, 
\begin{equation}\label{solution-pair}
  u^{\leq N}_\varepsilon(\theta)=\sum_{k=0}^N u_k(\theta)\varepsilon^k\quad
  \mbox{and}\quad
  c^{\leq N}(\varepsilon) =\sum_{k=0}^N c_k\varepsilon^k,
\end{equation}
such that $$\|E_{c^{\leq N}(\varepsilon)}[ u^{\leq N}_\varepsilon]\| = O(\e^{N+1})$$ where
$E$ is the functional defined in \eqref{func-eqn}.

\end{proposition}


\subsection{Pad\'e extrapolation}\label{pade}

The domain of analyticity for the solution of \eqref{func-eqn} can be
approximated by implementing a Pad\'e method
in which we use the approximate solutions obtained by the Lindstedt series
constructed in Section \ref{lindstedt}.

The Pad\'e method is quite standard and is presented in several places
in the literature.
Here, we follow the exposition in \cite{baker-96}. 
A Pad\'e approximant of order $[p/q]$ of a function $g(\e)=\sum_{i=0}^\infty g_i \e^i$ is a
rational function, $P(\e)/Q(\e)$, 
which agrees with $u$ to the highest possible order in $\e$.

That is, 
\begin{equation} \label{pade}
 g(\e)-\frac{P(\e)}{Q(\e)} = \mathcal{O}(\e^{p+q+1}).
\end{equation}

Where $P(\e)$ and $Q(\e)$ are polynomials of degree $p$
and $q$ respectively, $Q(0) = 1$.

The existence of the polynomials $P$ and $Q$ can be obtained by
noticing that \eqref{pade} is equivalent to
$$g(\e)Q(\e) = P(\e) + \mathcal{O}(\e^{p+q+1})$$
and, then, considering $P(\e)=\sum_{i=0}^p P_i\e^i$ and
$Q(\e)=\sum_{i=0}^q Q_i\e^i$ the coefficients of the polynomials
can be found by solving the following systems of equations 
\begin{eqnarray}\label{pade-linear}
u_i + \sum_{j=1}^iu_{i-j}Q_j =P_i\qquad 0\leq i \leq p \nonumber \\
u_i + \sum_{j=1}^qu_{i-j}Q_j =0 \qquad p < i \leq p+q.
\end{eqnarray}
The second equation of \eqref{pade-linear} gives
the $Q_j's$, and then we can find the $P_j's$ by substituting in the
first equation. Then, the boundary of the domain of
analyticity of a function can be approximated by the zeros of
$Q$ in the $[p/q]$ Pad\'e approximant.

There are a number of implementations of the Pad\'e methods that are used
in a quite standard manner. In the present work we use the implementations
included in Version 2.9.0 of GP/PARI, \cite{PARI98}.


\subsection{Newton's method}\label{newton}

In this section we summarize an iterative scheme in scales of Banach spaces that
can be very well adapted to perform numerical computations. The scheme is based in
a Newton iteration starting from approximate solutions
to the equation \eqref{func-eqn}. We briefly describe the scheme here since details
of schemes of these kind and numerical implementations have been described
already in the literature
\cite{Cal-Cel-10, Cal-Cel-Lla-16, Cal-Fig-12,Cal-Cel-Lla-13b}, and the reader
can refer to these works for more details.

We start from an approximate solution $(u_a, c_a)$ of equation \eqref{func-eqn}. Namely,
we have a solution so that $\|E_{c_a}[u_a]\|$ is small enough. The approximate solution
could be obtained by several means. One possibility is starting from the integrable
case (for $\e$ close to zero) and perfoming continuation
or from a Lindstedt series expansions like the
ones obtained in Section \ref{lindstedt}. We remark that in the dissipative
standard map we are studying, $\e =0$ is the point where the map becomes symplectic.
Since we use methods for conformally symplectic systems we actually start the
continuation from values of $\e$ that are note equal to zero but small.

The Newton algorithm consists of adding a correction $(v, \delta)$ to the approximate solution
so that supremum norm of \eqref{func-eqn} evaluated in the function plus the corrections,
$\|E_{c_a+\delta}[u_a+v]\|$, is of the order of the square of
the norm of \eqref{func-eqn} evaluated at the approximate solution,
\[\|E_{c_a+\delta}[u_a+v]\| \leq C \|E_{c_a}[u_a]\|^2.\]
One obtains the correction by solving the linearized equation of $E_{c_a+\delta}[u_a+v]$
for $(v, \delta)$  around the
approximate solution, $(u_a, c_a)$.

In this case, the equation we have to solve is
\begin{equation}
  D_u E_{c_a}[u_a] v - \delta = -E_{c_a}[u_a]
\end{equation}
which involves unbounded operators in Banach spaces (namely $D_u E_{c_0}[u_0] v$)
that are actually bounded if one considers that the operators map into Banach
spaces of less regularity.
It is a standard observation in Nash--Moser theory \cite{Zeh-75a, Zeh-76},
that to set up a converging iterative Newton scheme it
is not necessary to find an exact inverse of the operator
$D_u E_{c_a}[u_a]$, but an approximate inverse will suffice.

One obtains an approximate inverse by noticing that the modified
Newton equation,
\begin{equation}\label{eqn:newton-modified}
h'D_u E_{c_a}[u_a] v - v D_u E_{c_a}[u_a] h' = -h'(E_{c_0}[u_0]-\delta)\ ,
\end{equation}
with $h'(\theta) = 1 + \frac{\partial u(\theta)}{\partial \theta}$,
factorizes in a sequence of operators
that map Banach spaces of regular functions to Banach
spaces of functions with less regularity.

This method has been used in several works \cite{CellettiC09, Lla-08, Lla-01}.
Here we only make a reference to the justification in \cite{Cal-Cel-10}, where
the reader can refer to for details.
Let the operators $ \mathcal D_{-1}$, $\mathcal D_{1}^{b}$ by
\begin{eqnarray}\label{diff-eqns}
\mathcal D_{-} f(\theta) &=& f(\theta - \omega) - f(\theta)\nonumber\\
\mathcal D_{+}^b f(\theta) &=& f(\theta + \omega) - b f(\theta)\ .
\end{eqnarray}
A small remark is that \ref{diff-eqns} are operators that are diagonal
in Fourier space.
In the following lemma, we write the modified Newton as
a sequence of operators
that are either diagonal in Real or Fourier space. 

\lem{factorization}
The modified Newton equation in \eqref{eqn:newton-modified}
with $E_{c_a}[u]$ defined in \equ{func-eqn} is equivalent to
\begin{equation}\label{equ:qp-reduced-newton}
\mathcal D_+^b[-h'(\theta)h'(\theta - \omega) \mathcal D_{-} [(h')^{-1}(\theta) v(\theta)
] ] = - h'(E_{c_a}[u_a](\theta) - \delta).
\end{equation}

\elem

\begin{remark}\label{diagonal}
One notices that the operators involved in the l.h.s. of equation
\eqref{equ:qp-reduced-newton} only involve differentiation,
multiplication, division, shifting the arguments of functions, and
solving the difference equations with constant coefficients in \eqref{diff-eqns}.
All this operations can be implemented very efficiently using the computer.
For instance if we discretize the periodic functions using $n$ uniformly distributed
points and we use a Fast Fourier Transform method, the modified Newton step
equation can actually be solved in $O(n\log n)$ operations. 
\end{remark}

The factorization in equation
\eqref{equ:qp-reduced-newton} suggests an algorithm that is
used to solve the modified Newton equation.  

\begin{algorithm}\label{algorithm2}
  
  \begin{enumerate}

  \item[\textbf{i)}]
    Find two functions $\varphi$ and $\nu$ solving the equations
\begin{equation}
   \mathcal D^b_{+}\varphi (\theta) =  -h'E_{c_a}[u]
   \label{firstsolution1}
\end{equation}
and
\begin{equation}
   \mathcal D^b_{+}\nu (\theta) =  -h'(\theta)\ .
   \label{firstsolution2}
\end{equation}
Notice that if $\varphi(\theta)$  and $\nu(\theta)$ are solutions of
\eqref{firstsolution1} and \eqref{firstsolution2}, respectively, then
the equation $\mathcal D^b_+(\varphi(\theta) - \delta \nu(\theta)) =
-h'(\theta) (E_{{c_a}}[u_0](\theta)-\delta)$
holds for any $\delta \in \mathbb{C}$. This will allow us to chose a
complex number $\delta$ so that the average of
$\frac{\varphi(\theta) - \delta \nu(\theta) }{h'(\theta)h'(\theta - \omega)}$
vanishes.

\item[\textbf{ii)}]
  Choose $\delta \in \mathbb{C}$ such that \[\avg{\frac{\varphi(\theta) -
      \delta \nu(\theta) }{h'(\theta)h'(\theta - \omega)}} = 0\ .\]

\item[\textbf{iii)}]
  Obtain $w$ from the solution of the constant coefficient difference
equation
\begin{equation}\label{CpGzeroaverage2}
  \mathcal D_{-} w (\theta) = \frac{\varphi(\theta) -
    \delta \nu(\theta)}{-h'(\theta)h'(\theta - \omega)}\ .
\end{equation}
Notice that after choosing a $\delta$ in step ${\bf ii)}$ so that the right hand side has zero
average we can always find a periodic function
$w$ solving \eqref{CpGzeroaverage2}
when the r.h.s. is smooth enough.

\item[\textbf{iv)}]
  Construct $v(\theta) = h'(\theta)w(\theta)$ and obtain the improved solution
$(\tilde u, \tilde c)$ defined as
$$
\tilde{u}(\theta) = u_a(\theta) + v(\theta)\ ,\qquad \tilde{c} = c_a + \delta\ .
$$

\end{enumerate}

\end{algorithm}
The observation in remark \ref{diagonal} is that the operators in
\eqref{equ:qp-reduced-newton} are very efficiently implementable
with the use of a computer either in Real or in Fourier space.
This efficiency comes from the fact that all the
operations involved in the
four steps of
Algorithm \ref{algorithm2} are multiplications, additions and integrals
of periodic functions that take only $O(n)$ operations in Real
space; and differentiation, shifts and solving cohomology equations
with constant coefficients, that take only $O(n)$ operations in Fourier
space. Therefore, the most expensive operation in the Algorithm
\ref{algorithm2} is transforming from Real to Fourier
space and back. This can be done in $O(n \log n)$ operations by means of
a Fast Fourier Transform. 

\begin{remark}
  We note that the algorithm is guaranteed to converge inside the boundaries of the
  analyticity domain. Indeed, in \cite{Cal-Lla-10} was rigorously justified
  that the algorithm only fails to converge as the continuation reaches the boundary
  of analyticity. Therefore, the continuation method can also be used to asses the
  bounds on the domain of $\e$.
\end{remark}

\section{Numerical results}

In this section we present the results of implementing the methods described
in Section \ref{methods}. All the computations were done using the golden ratio, $\omega= \frac{1+\sqrt{5}}{2}$, which satisfies \eqref{dioph-cond} \cite{Lla-01}.

\subsection{Lindstedt expansions}
The construction of Lindstedt series in Section \ref{lindstedt} was implemented
as a numerical algorithm.
The statement of Proposition \ref{prop:existenceLind}
tells us that given any $N\in \mathbb N$, the outcome of the method is the pair of polynomials
of degree $N$  in \eqref{prop:existenceLind}.
The observation of Lemma \ref{L_w_sol} is that the operator
$L_\omega$ defined in equation \eqref{Lomega}
is diagonal in Fourier
Space and equation \ref{solveLomega} can be solved for $\phi$ if we allow
to obtain functions with less regularity than the rigth hand side, $\eta$.
We find the solution numerically by transforming
to Fourier space and solving for the $u_k$'s from expressions
\eqref{ec_u1} to \eqref{ec-u_n}. At every order of the process we obtain
the $c_k$'s as a byproduct of imposing
the condition that every order should have zero average.

The Lindstedt series expansions are used to obtain an approximate solution
to the functional equation in
\eqref{func-eqn} at some high order. Indeed, we discovered that with our implementations
it is very hard to notice that the functions are not analytic. Namely, the singularities
that exist close to the point $\e = 0$ are very hard to detect so we have no evidence that
the radius of convergence of the series is exactly zero in the complex plane.
Thus, if the solution belongs to a Gevrey class then the Gevrey exponent would be very close to one.

We approximated several norms of the coefficients, $u_k(\theta)$,
to have an indication of how far the functions are from being analytic.
First, we use the
norm on the complex strip of size $\rho > 0$, i.e., 
$\theta \in \mathbb S^1_\rho$ if $|\Im (\theta)|<\rho$.
Let $f: \mathbb S^1_\rho \to \mathbb S^1_\rho$ be a function of $\mathbb S^1_\rho$
then the norm we use is  
\[\|f\|_\rho = \sum_{\ell\in \mathbb Z} |\hat f_\ell |^2 e^{2 \pi |\ell| \rho}\]
where $\hat f_\ell$ are the Fourier coefficients of $f$.

We say that the function $f(\e)$ belongs to the Gevrey class $G^\sigma$ with respect to
the norm $\| \cdot\|_B$ at $\e=0$ whenever
\[ \|\partial_\e^k f(\e)\|_B \leq C R^k k^{\sigma k},\]
for $\e = 0$, \cite{Gevrey1918}.


Since we want to check if the function $u_\e(\theta)$ belongs to a Gevrey
class at $\e=0$ with the analytic norms it is convenient to
compute the following expressions as
functions of $k$,
\begin{equation}\label{analytic-gevrey}
  A_\rho(k) \equiv \frac{1}{k} \log  \|u_k(\theta)\|_\rho,
\end{equation}
and then approximate the constant $\sigma$.

\begin{figure}[h]
\includegraphics[width=15truecm]{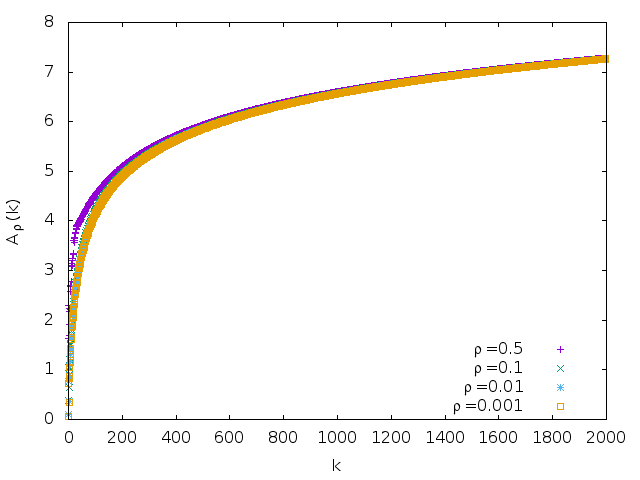}
\caption{Analytic norms of the coefficients of the Lindstedt expansion
  plotted as in expression \eqref{analytic-gevrey}.}
\label{analytic-norms}
\end{figure}

The expressions \eqref{analytic-gevrey} as functions of
$k$ for the coefficients of the approximate solution are shown in
Figure \ref{analytic-norms}.

We also used Sobolev norms defined for a real number $r>0$ by
the $L^2$-norm of the $r^{th}$ derivative with respect
to $\theta$,
\[ \|f\|_r = \|\partial^r_\theta f\|_{L^2}.\]
Notice that when $r=0$, $\|f\|_0$ corresponds to the $L^2$ norm of $u$.
The Sobolev norms can also be written in terms of Fourier coefficients
as follows,
\[\| f \|_r = \left ( \sum_{k\in \mathbb Z} (2 \pi k)^{2r} |\hat f_k|^2\right)^{1/2}.\]

\begin{figure}[h]
\includegraphics[width=15truecm]{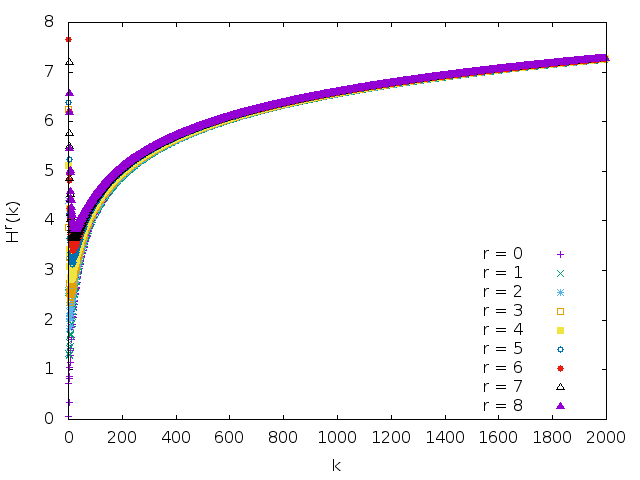}
\caption{Sobolev norms of the coefficients of the Lindstedt expansion
  plotted as functions of the order of
  $\e$.}
\label{sobolev-norms}
\end{figure}

As in the case for analytic norms we define the following expressions for the
Sobolev norms,
\begin{equation}\label{sobolev-gevrey}
  H^r(k) \equiv \frac{1}{k} \log  \|u_k(\theta)\|_r.
\end{equation}
We include the values of $H^r(k)$ for the coefficients of the approximate solution
and several values of $r$
in Figure \ref{sobolev-norms}.

In both cases, the behavior of the norms coefficients $\|u_k(\theta)\|_B$ with respect
to $k$ seem to belong to Gevrey classes. In Tables I and II we include the fit of the
plots in Figures \ref{analytic-norms} and \ref{sobolev-norms}.

\begin{table}\label{tab-analytic}
\begin{tabular}{ |p{2cm}||p{4cm}|p{4cm}|p{4cm}|  }
 \hline
 \multicolumn{4}{|c|}{$A_\rho(k) = log(a) + c \log(k+b)$} \\
 \hline
 & a & b & c\\
 \hline
 $\rho = 0.5$   &  0.719467892188978    & 27.3937734029272 & 1.00001766652951\\
 $\rho = 0.1$ &   0.719927523693294  & -6.67612303717059 & 0.999991308925563\\
 $\rho = 0.01$ & 0.719826182759749  &  -12.4920865977342 &  0.999998383397402 \\
 $\rho = 0.001$ & 0.719813322136991 & -13.0568139120788 & 0.999999331056825 \\
 \hline
\end{tabular}
\caption{Numerical fit of analytic norms in expression \eqref{analytic-gevrey}
for different values of $\rho$.}
\end{table}


\begin{table}
\begin{tabular}{ |p{2cm}||p{4cm}|p{4cm}|p{4cm}|  }
 \hline
 \multicolumn{4}{|c|}{$H^r(k) = log(a) + c \log(k+b)$} \\
 \hline
    & a & b & c\\
 \hline
 $r = 0$ &  0.71981186150290  & -13.11936925724 & 0.99999943903933 \\
 $r = 1$ &  0.71990246696872  & -8.455098498275 & 0.99999292602721 \\
 $r = 2$ &  0.71995774061454  & -3.699094566292 & 0.99998926987818 \\
 $r = 3$ &  0.71997659277593  &  1.151414934192 & 0.99998826853277 \\
 $r = 4$ &  0.71995790042650  &  6.099307051805 & 0.99998972777892 \\
 $r = 5$ &  0.71990050891201  &  11.14758050573 & 0.99999346098676 \\
 $r = 6$ &  0.71980323117743  &  16.29938193670 & 0.99999928901585 \\
 \hline
\end{tabular}
\label{tab-sobolev}
\caption{Numerical fit of analytic norms in expression \eqref{sobolev-gevrey}
for different exponents, $r$.}
\end{table}

The numerical results in Tables I
and II lead us to think that the solutions that we approximate
are funcions that are very hard to distinguish from analytic functions
by just examining the a truncated expansion series.
One of the rigorous results in \cite{ Cal-Cel-Lla-16} states that
the functions that satisfy equation \eqref{func-eqn} fail to be analytic
since there is no ball around $\e = 0$ where the formal power series converges.
Therefore, we conjecture that the solutions belong to a Gevrey class with
an exponent that is very close to the analytic class.
\begin{conjecture}
The parametrization $u_\e$ belongs to a Gevrey class, $G^\sigma$, as a function of $\e$. The index, $\sigma$, is close to $1$. 
\end{conjecture}

\subsection{Approximation of poles of the Lindstedt series}\label{PadeApprox}

\begin{figure}[H]
\includegraphics[width=10truecm]{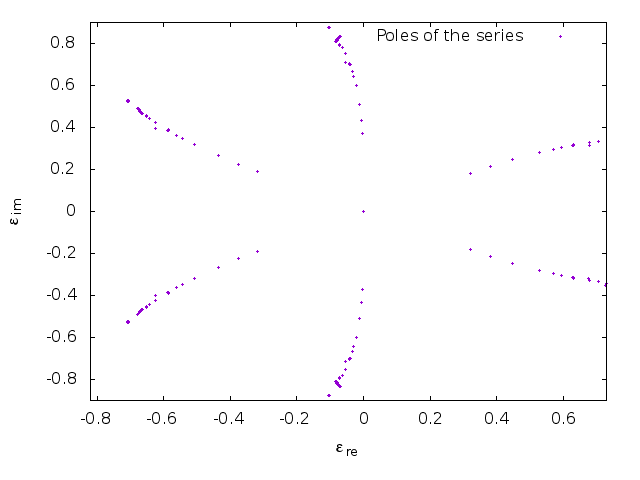}
\includegraphics[height=7.5truecm]{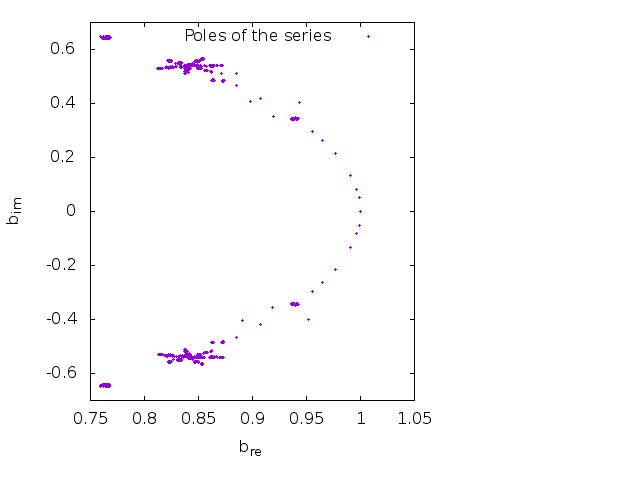}
\caption{Points which are simultaneously poles of Pad\'e approximants of degree [475,475] and [500,500]. The implementation was done with 1000 digits.
  Left panel: Poles in the complex plane $\e \in \mathbb C$.
  Right panel: Poles evaluated in the function $b(\e) = 1 - \e^3$, with  $\e \in \mathbb C$. }
\label{poles1}
\end{figure}

\begin{figure}[h]
  \includegraphics[height=8truecm]{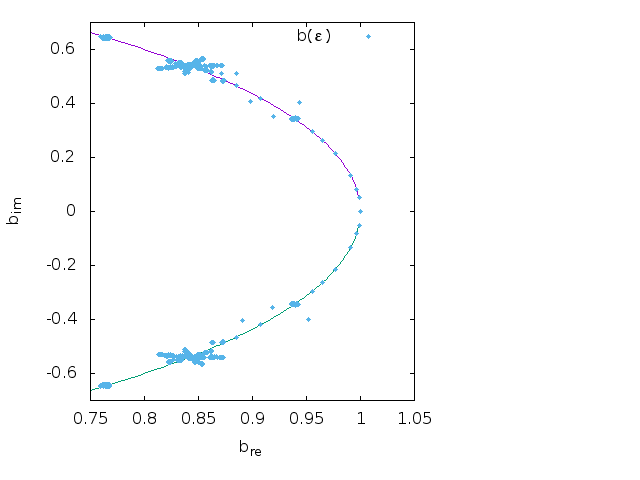}
  \includegraphics[width=7truecm]{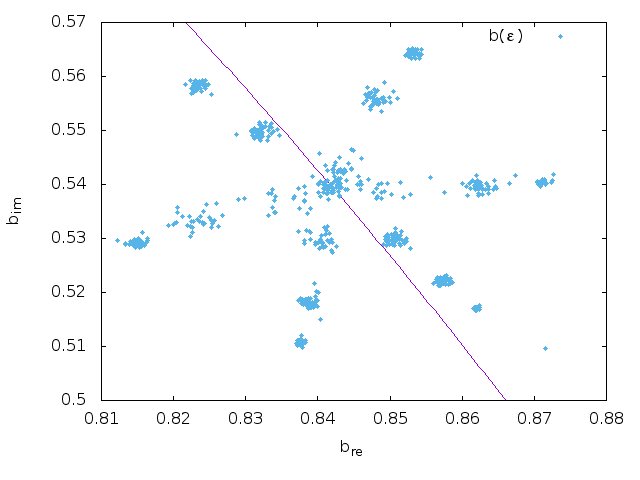}
  \includegraphics[width=7truecm]{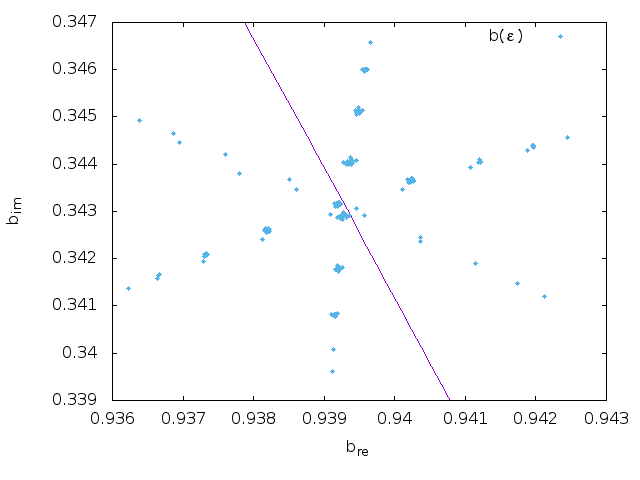}
\caption{The poles compared to the unit circle.
  Upper panel: Evaluation of the poles of the series by the function $b(\e) = 1- \e^3$.
  Lower panels: Two zoomed in versions of the set.}
\label{poles2}
\end{figure}

Here we include the poles of the Lindstedt polynomial
found with the Pad\'e method.
In Figure \ref{poles1}, we show the poles of the series approximated by means of
the Pad\'e method.
It is well known that the Pad\'e method
computations are very sensitive to precision,
see \cite{baker-96}, so
we have implemented the computations with extended precision
using the software gp/Pari \cite{PARI98}.
We show the values of the poles in the $\e$ complex plane,
and the complex values of the function $b(\e) = 1 - \e^3$.
Figure \ref{poles2}, contains the comparison of the values
of the function $b(\e)$ with the unit circle. We also include
zoomed in versions of the values of $b(\e)$ in Figure \ref{poles2}.


\subsection{Newton method}
We used Newton's method and continuation to explore the monodromy of the
solutions in the domains. A rigorous result in \cite{Cal-Cel-Lla-16} states
that the solutions defined in the domain of analyticity in $\e$ have
trivial monodromy. We verified this fact numerically
by perfoming continuation of the solutions $(u_\e, c_\e)$ around the poles
that were previously approximated using the Pad\'e series method
described in Section \ref{PadeApprox}.

We used the approximated poles as centers of circular paths in $\e$
over which we perfomed continuation while solving the invariance
equation \eqref{func-eqn} using Algorithm \ref{algorithm2}.
Once the continuation completes a complete turn around a chosen pole,
one verifies that the solution always arrives to the same starting point.
This is an effect of the monodromy of the functions being trivial.

\begin{figure}[h]
\centering
\includegraphics[height=10truecm]{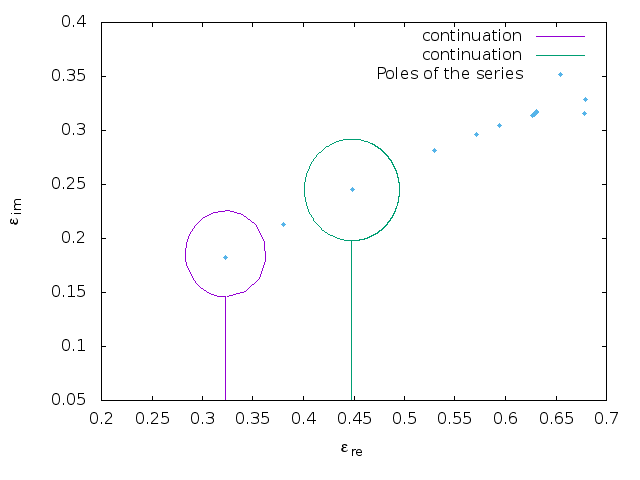}
\caption{Poles of the series and two different continuations done with
  the Newton algorithm. The continuation is done around the pole in order
  to illustrate that the monodromy is trivial. }
\label{no-mon-path}
\end{figure}

\begin{figure}[h]
\centering
\includegraphics[height=6truecm]{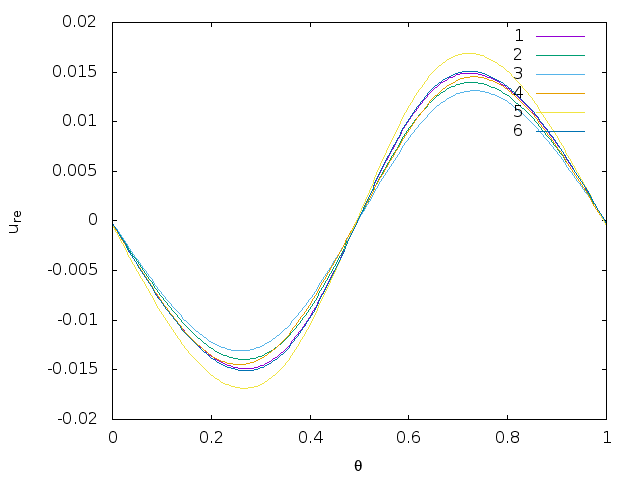}\quad
\includegraphics[height=6truecm]{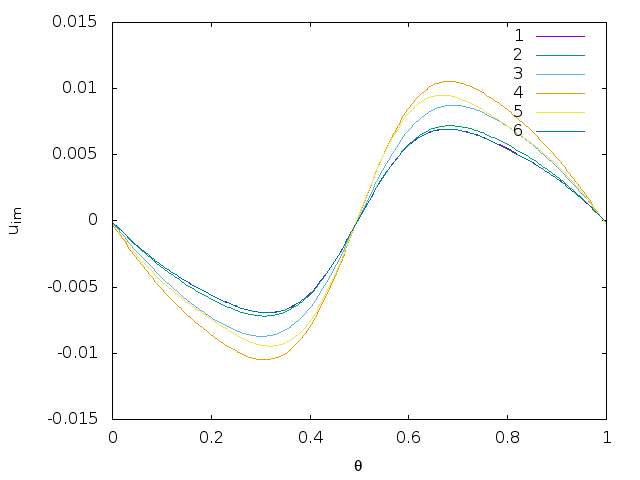}
\caption{Real and imaginary part of different instances of
  a continuation by the Newton algorithm including the initial and final functions.
  One observes that there is no monodromy after a full turn around the pole.}
\label{no-mon-re-im}
\end{figure}

We present several
instances of the functions for different parameter
values along a circle winding around a pole in Figure \ref{no-mon-re-im}.
The path we used to surround the pole is presented in Figure \ref{no-mon-path}.
The continuation was perfomed using FFTW3, \cite{FFTW05}, with the
libquadmath library, \cite{QuadMath}.
The radii of the continuation paths
were chosen so that the path did not come very close to the poles. Indeed,
when the continuation comes close to a pole our implementation of the
Newton method becomes degenerate in the sense that one needs to compute
quotients of very small quantities. The reason is that when solving
equations \eqref{firstsolution1} and \eqref{firstsolution2},
the divisors depending on $\e$, are below machine precision
close to the pole and
dividing over those quatities leads to large numerical errors.






\vspace{-10pt}

\begin{table}
\begin{tabular}{|c|c|c|}
 \hline
 \hline
 Instance & $\e$ & $c(\e)$\\
 \hline
 
1 & $0.3202966+i0.1460915$ & $0.01994937-i0.06774761$ \\

2 & $0.3008391+i0.1527000$ & $0.009976542-i0.06136120$ \\

3 & $0.2830167+i0.1871540$ & $-0.01146081-i0.06221038$ \\

4 & $0.3122423+i0.2245263$ & $-0.02718298-i0.08804174$ \\

5 & $0.3613448+i0.1973876$ & $0.007928831-i0.1127768$ \\

6 & $0.3242691+i0.1460201$ & $0.02157160-i0.06953580$ \\

 \hline
\end{tabular}
\label{tab-sobolev}
\caption{Values of $\e$ and $c(\e)$ for different instances taken from the small circle in FIG.\ref{no-mon-path}.}
\end{table}

\newpage
\noindent
{\bf Acknowledgments.} We thank Rafael de la Llave and Alessandra Celletti 
for fruitful discussions.
This work was started while A.B. was a student at IIMAS and was supported
by CONACYT.

\bibliographystyle{alpha}
\bibliography{DomComp}

\end{document}